\documentclass[11pt]{amsart}
\usepackage{epsfig,amsmath,amssymb,amsthm,epsf}
\usepackage{graphicx,array}

\newcommand{\zz}[1]{\mathbb #1}
\newtheorem{proposition}{Proposition}
\newtheorem{lemma}{Lemma}

\newtheorem{theorem}{Theorem}
\newtheorem{claim}{Claim}

\begin{document}
\title{ An oriented competition model on $Z_{+}^2$. }
\author{{George Kordzakhia}}
\address{University of California\\
Department of Statistics\\
Berkeley CA }
\email{kordzakh@stat.berkeley.edu}
\author{Steven P. Lalley} 
\address{University of Chicago\\ Department of Statistics \\ 5734
University Avenue \\
Chicago IL 60637}
\email{lalley@galton.uchicago.edu}
\date{\today}
%\maketitle
\begin{abstract}
We consider a two-type oriented competition model on the 
first quadrant of the two-dimensional integer lattice. 
Each vertex of the space may contain only one particle of either Red
type or Blue type.  A vertex flips to the color of a randomly chosen
southwest nearest neighbor at exponential rate $2$.
At time zero there is one Red particle located at $(1,0)$ and one Blue
particle located at $(0,1)$.  The main result is a partial shape
theorem: Denote by $R (t)$ and $B (t)$ the red and blue regions at
time~$t$. Then (i) eventually the upper half of the unit square
contains no points of $B (t)/t$, and the lower half no points
of $R (t)/t$; and (ii) with positive probability there are angular
sectors rooted at $(1,1)$ that are eventually either red or blue.  The
second result is contingent on the uniform curvature of the boundary
of the corresponding Richardson shape.
\end{abstract}
\maketitle
Key words: competition, shape theorem, first passage percolation.

\section{Introduction.}\label{sec:intro} In this paper we study a
model where two species Red and Blue compete for space on the first
quadrant of $\mathbb{Z}^2$.  At time $t>0$ every vertex of $\zz{Z}^{2}$
is in one of the three possible states: vacant, occupied by a Red
particle, or occupied by a Blue particle.   An unoccupied vertex $z= (x,y)$ may be
colonized from either $(x,y-1)$ or $(x-1,y)$ at rate equal to the
number of occupied south-west\ neighbors; at the instant of first
colonization, the vertex flips to the color of a randomly chosen
occupied south-west neighbor. Once occupied, a vertex remains
occupied forever, but its color may flip: the flip rate is equal to
the number of south-west neighbors occupied by particles of the
opposite color.  The state of the system at any time $t$ is given by
the pair $R(t),B(t)$ where $R(t)$ and $B(t)$ denote the set of sites
occupied by Red and Blue particles respectively. The set $R (t)\cup B
(t)$ evolves precisely as the occupied set in the oriented Richardson
model, and thus, for any initial configuration with only finitely many
occupied sites, the growth of this set is governed by the \emph{Shape
Theorem}, which states that the set of occupied vertices scaled by time
converges to a deterministic set $\mathcal{S}$ (see for example
\cite{cox}). A rigorous construction and more detailed description of
the oriented competition model is given in Section \ref{sec:Richardson}.

\begin{figure}
\label{OrientedPicture} \vspace{1.5ex}
\begin{tabular}{cc}
\includegraphics[width=2.5in, height=2.5in]{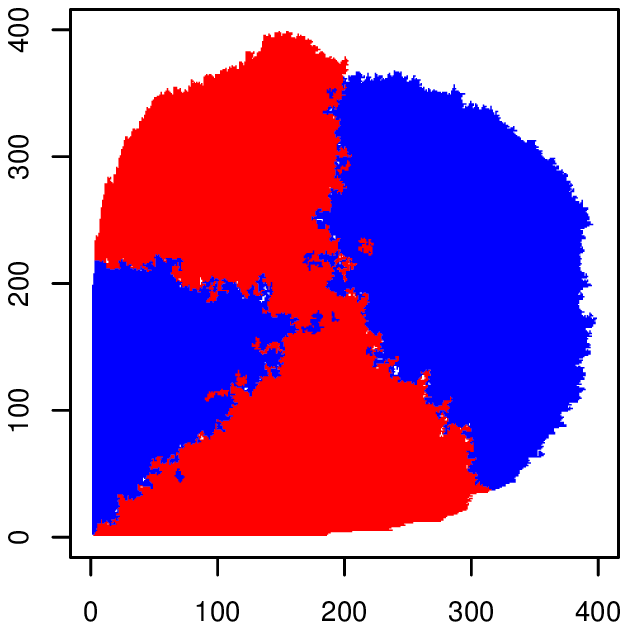}
\includegraphics[width=2.5in, height=2.5in]{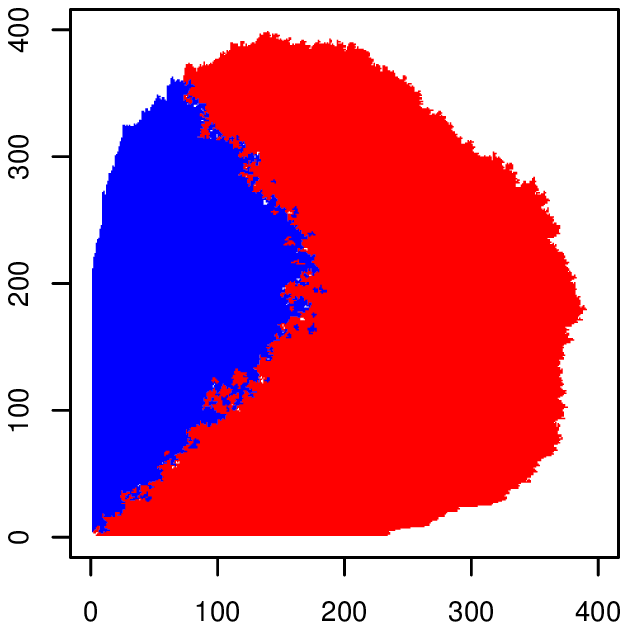} 
\end{tabular}
\caption{Two Realizations of  Oriented Competition.}
\end{figure}

The simplest interesting initial configuration has a single Red
particle at the vertex $(1,0)$, a single Blue particle at $(0,1)$, and
all other sites unoccupied. We shall refer to this as the
\emph{default} initial configuration. When the oriented competition
process is started in the default initial configuration, the red and
blue particles at $(1,0)$ and $(0,1)$ are protected: their colors
can never be flipped. Thus, both colors survive forever w.p.1.
Computer simulations for the oriented competition model started in the
default and other finite initial configurations suggest that the
shapes of the regions occupied by the Red and Blue types stabilize as
times goes to infinity -- see Figure \ref{OrientedPicture} for snapshots of
two different realizations of the model, each started from the default
initial configuration. A peculiar feature of the
stablization is that the limit shapes of the red and blue regions are
partly deterministic and partly random: The southeast corner of the
occupied region is always equally divided between the red and blue
populations, with boundary lying along the line $y=x$. However, the
outer section seems to stablize in a random union of angular wedges
rooted at a point near the center of the Richardson shape. Although
the location of the root appears to be deterministic, both the number
and angles of the outer red and blue regions vary quite dramatically
from one simulation to the  next. 

The purpose of this paper is to prove that stabilization of Red
and Blue zones occurs with positive probability. (We conjecture that in
fact it occurs with probability $1$, but we have been unable to prove
this.) To state our result precisely, we shall need several facts
about the limit shape $\mathcal{S}$ of the oriented Richardson model
$Z (t):= R (t)\cup B (t)$. The proof of the Shape Theorem \cite{cox}
shows that $\mathcal{S}$ is a compact, convex subset of the first
quadrant of $\zz{R}^{2}$. It is generally believed -- but has not been
proved -- that the outer boundary $\partial^{o}\mathcal{S}$ of
$\mathcal{S}$ (the portion of $\partial \mathcal{S}$ that lies in the
interior of the first quadrant) is \emph{uniformly curved}, that is,
for every point $x$ in this part of the boundary there is a circle of
finite radius passing through $x$ that contains $\mathcal{S}$ in its
interior. In section \ref{sec:environment} we shall prove the following.

\begin{lemma}  \label{ProperSubset} 
The Richardson shape $\mathcal{S}$ has the points $(1,0)$ and $(0,1)$
on its boundary, and the point $(1,1)$ in its interior.
\end{lemma}

\noindent It will follow by convexity that the unit square
$\mathcal{Q}=[0,1]^{2}$ lies entirely in $\mathcal{S}$. Define
$\mathcal{Q}_{1}$ and $\mathcal{Q}_{2}$ to be the subsets of
$\mathcal{Q}$ that lie (strictly) above and below the main diagonal
$x=y$. \\
\indent  For any subset $Z \subset   \mathbb{R}^2$, define 
\begin{equation*}
\hat{Z}=\{x \in \mathbb{R}^2: \mbox{ dist}(x,Z) \le 1/2 \},
\end{equation*}
where dist denotes distance in the $L^{\infty}$-norm on
$\mathbb{R}^2$.  For any set $Z \subset \mathbb{R}^2$ and any scalar
$s>0$, let $Z/s= \{ y/s: y \in Z \}$. 

%%%%%%%%%%%%%%%%%%THEOREM
\begin{theorem}  \label{MainTheorem}
With probability one, for all large $t$
\begin{equation}\label{eq:q1q2} 
\mathcal{Q}_1   \subset \hat{R}(t) /t \quad \text{and} \quad 
 \mathcal{Q}_2 \subset \hat{B}(t) /t.  
\end{equation}
Furthermore, if the outer boundary $\partial^{o}\mathcal{S}$ of the
Richardson limit shape $\mathcal{S}$ is uniformly curved, then for
every $\epsilon>0$  the following holds with positive
probability: There exist random angular sectors $A_1,..,A_n$ rooted at
$(1,1)$ that do not intersect the open unit square $\mathcal{Q}^{o}$
such that 
\begin{enumerate}
\item [(a)] eventually $A_i$ is either Red or Blue, and 
\item [(b)]  the complement of $\bigcup A_i$ in $\mathcal{S}
\setminus \mathcal{Q}$ has angular measure less than $\epsilon$.
\end{enumerate}
\end{theorem}
%%%%%%%%%%%%%%%%%%%%%%
\indent Another competition model on $\mathbb{Z}^d$ (non-oriented
version) was studied in~\cite{kordzakh}. It was shown that if the
process starts with finitely many particles of both types (Red and
Blue), then the two types coexist with positive probability under the
condition that the shape set of the corresponding non-oriented
Richardson model is uniformly curved. The behavior of the oriented
model differs from that of the model considered in \cite{kordzakh} in
that the limit shape contains the determinsitic component (\ref{eq:q1q2}).
%%%%%%%%%%%%%%%%%%%%%%%%%%%%%%%%%%%%%%%%%%%%%%%%%%%%%%%%%%%%%%%%%%%%%%%%%%%%%%%%%%%%%%%%%%%%%%%%
\section{Preliminaries}\label{sec:preliminaries}
\subsection{Graphical Constructions}\label{sec:graphicalConstructions}
%%%%%%%%%%%%%%%%%%%%%%%%%%%%%%%%%%%%%%%%%%%%%%%%%%%%%%%%%%%%%%%%%%%%%%%%%%%%%%%%%%%%%%%%%%%%%%%
The competition model, the Richardson model, and the competition model
in a hostile environment may be  built using the same  \emph{percolation structure}  $\Pi$.
For details on percolation structures  see \cite{DurrettBook}. 
Here we briefly describe the construction of $\Pi$.
To each directed edge $xy$ 
%of the lattice $\zz{Z}_{+}^{2}$
such that $x \in \zz{Z}_{+}^{2} \setminus \{(0,0)\} $, and
$y=x+(0,1)$ or $y=x+(1,0)$ is assigned   rate-$1$ Poisson process. 
The Poisson processes are mutually independent.
 Above each vertex $x$ is drawn a timeline,
 on which are placed marks at the
occurrence times $T^{xy}_{i}$ of the Poisson processes attached to
directed edges emanating from $x$; at each such mark, an arrow is
drawn from $x$ to $y$.
 A \emph{directed path} through the percolation
structure $\Pi$ may travel upward, at speed $1$, along any timeline,
and may (but does not have to) jump across any outward-pointing arrow
that it encounters. A \emph{reverse path} is a directed path run
backward in time: thus, it moves downward along timelines and jumps
across inward-pointing arrows. A \emph{voter-admissible} path is a
directed path that does not pass any inward-pointing arrows. Observe
that for each vertex $z$ and each time $t>0$ there is a unique
voter-admissible  path beginning at time $0$ and terminating
at $(z,t)$: its reverse path is gotten by traveling
downward along timelines, starting at $(z,t)$, jumping across all
inward-pointing arrows encountered along the way.  \newline
\indent For each $(z,t)$ denote by $\Gamma (z,t)$ the collection of reverse
paths on  percolation structure $\Pi$  originating at $(z,t)$ 
and terminating in $(\mathbb{Z}_{+}^2,0)$.  We also use 
$\Gamma (z,t)$ to denote the set of ends of all paths in the collection.
There exists a unique reverse voter-admissible path 
$\tilde{\gamma}_{(z,t)}=\tilde{\gamma}$ in $\Gamma (z,t)$. 
We say that a path $\gamma$ has {\it attached end},
or $\gamma$  is {\it attached}, if it terminates in $(R(0) \cup B(0),0)$.
For $s \in [0,t]$  denote by $\gamma(s)$ the location of the path
in $\mathbb{Z}_{+}^2$ at time $t-s$, 
i.e. $\gamma(s)=z'$ if $(z', t-s) \in \gamma$. 
We can  now put an order relation on $\Gamma (z,t)$ as follows.
For two reverse paths $\gamma_1, \gamma_2 \in  \Gamma (z,t)$, 
let $\tau_{i}=\inf \{ s>0: \gamma_{i}(s)=\tilde{\gamma}(s) \},\ i=1,2$,
and set  $\gamma_1 \prec \gamma_2$ if $\tau_1 \le \tau_2$.
The order relation sets a priority on assigning an { \it ancestor}.
 A vertex $z$ is occupied by a particle at time $t$ if and only if
there is at least one attached reverse path originating at $(z,t)$. 
The set of terminating points  of all attached  paths in $\Gamma (z,t)$ 
is referred as the {\it set of potential ancestors} of the particle at $(z,t)$.  
The maximal element $\hat{ \gamma}$ in the set of attached  paths uniquely
determines the ancestor. Let $(z',0)$ be the terminating point  of
$\hat{ \gamma}$. Then the particle at $(z',0)$ is said to be the {\it ancestor} 
of the particle at $(z,t)$. 
Note that $(z,t)$ is vacant if and only if the set of attached paths is empty.
%%%%%%%%%%%%%%%%%%%%%%%%%%%%%%%%%%%%%%%%%%%%%%%%%%%%%%%%%%%%%%%%%%%%%%%%%%%%%%%%%%%%%%
\subsection{The simplest oriented growth model.}\label{sec:Richardson}
%%%%%%%%%%%%%%%%%%%%%%%%%%%%%%%%%%%%%%%%%%%%%%%%%%%%%%%%%%%%%%%%%%%%%%%%%%%%%%%%%%%
Denote by $Z(t)$ the set of vertices occupied  by time $t$,
and fix an initial configuration $Z(0)=\{ (0,1),(1,0) \}$ .
The Richardson model can be built using percolation structure as follows.
Set $Z(t)$ to be the set of vertices $z$ in $\mathbb{Z}^2_{+}$ such
that there is a directed path in  $\Pi$  that 
starts at $(Z(0),0)$ and terminates at $(z,t)$. 
% Denote by $T(z)$ the time at which $z$ gets occupied  
For $z \in \mathbb{R}^2_{+}$ let
$T(z)=\inf \{t: z \in \hat{Z}(t) \}$.
% (if $z$ has non-integer coordinates the occupation time is set 
% to  be the occupation time of the nearest  element of  $\mathbb{Z}^{2}_{+}$), 
and let $ \mu(z)= \lim_{n \rightarrow \infty} n^{-1}T(nz)$. 
The limit exists almost surely  by subadditivity.   
The growth  of $\hat{Z}(t)$ is governed by a Shape Theorem. 
A weakened version of the standard Shape theorem  may be obtained 
by using  subadditivity arguments.
The problem with a standard version of the Shape  theorem was 
that  $\mu(z)$ was not known to be
continuous on the boundaries of  $\mathbb{R}^{2}_{+}$.
In  \cite{martin}  J. Martin showed  that   $\mu(z)$ is
continuous on all  of  $\mathbb{R}^{2}_{+}$, and established the Shape theorem.
Furthermore, large deviations results for the Richardson model follow from 
papers by  Kesten \cite{kesten} and Alexander\cite{alexander}.
\begin{theorem} \label{OrientedShapeTheorem}
There exists  a non-random compact convex subset $\mathcal{S}$ of $\mathbb{R}_{+}^2$
such that for $\alpha \in (1/2,1)$, constants $c_1,c_2>0$ (depending on $\alpha$)
and  all $t>0$
\[
 P[ \mathcal{S}(t-t^{\alpha})  \subset \hat{Z}(t)   \subset \mathcal{S}(t+t^{\alpha})]>
  1-c_1 t^2  \exp\{ -c_2 t^{(\alpha-1/2)} \}.
\]
\end{theorem}
% Note that $\mu$ is a norm on  $\mathbb{R}^{2}_{+}$, and
% $\mu(z)= inf\{ t : z \in \mathcal{S}t \}$.
Let $\tilde{ \mathcal{S} }$ be  the limit set of the South-West 
oriented Richardson model.
This process starts with two particles at the vertices 
$(-1,0)$ and $(0,-1)$, 
and lives in the third quadrant of $\mathbb{Z}^2$.
It is easy to see that  $\tilde{\mathcal{S}}= - \mathcal{S}$.
For $\epsilon>0$ define the cone $K_{\epsilon}$ rooted  at $(1,1)$ 
by 
\[ 
K_{\epsilon}=\{ z \in \mathbb{R}^2: \mbox{ arg}\{z-(1,1) \} \in 
(-\pi/2+\epsilon,  \pi-\epsilon)\}.
\]
The following lemma follows from an elementary geometric argument.
The proof is identical to the proof of  Lemma 4 in \cite{Newman}. 
\begin{lemma} \label{OrientedShapeIntersection}
Suppose that 
 $\partial^{o}\mathcal{S}$  is uniformly curved.
For every  $\epsilon>0$ and  $\alpha \in (1/2,1)$ 
there exists  $c>0$ so that 
if  $z \in \partial{ \mathcal{S}} \cap K_{\epsilon}$,  then
for all $t_1,t_2>0$  we have
\[
\mathcal{S}(t_1+t_1^{\alpha}) \cap 
(z(t_1+t_2)+   \tilde{\mathcal{S}} (t_2+t_2^{\alpha})     )
 \subset D(zt_1,c(t_1+t_2)^{(\alpha+1)/2}).
\]
\end{lemma}
%%%%%%%%%%%%%%%%%%%%%%%%%%%%%%%%%%%%%%%%%%%%%%%%%%%%%%%%%%%%%%%%%%%%%%%%%%%%%%%%%
\section{Growth and competition in hostile environment.} \label{sec:environment}
\begin{figure}
\label{Hostile} \vspace{1.5ex}
\begin{tabular}{cc}
\includegraphics[width=2.5in, height=2.5in]{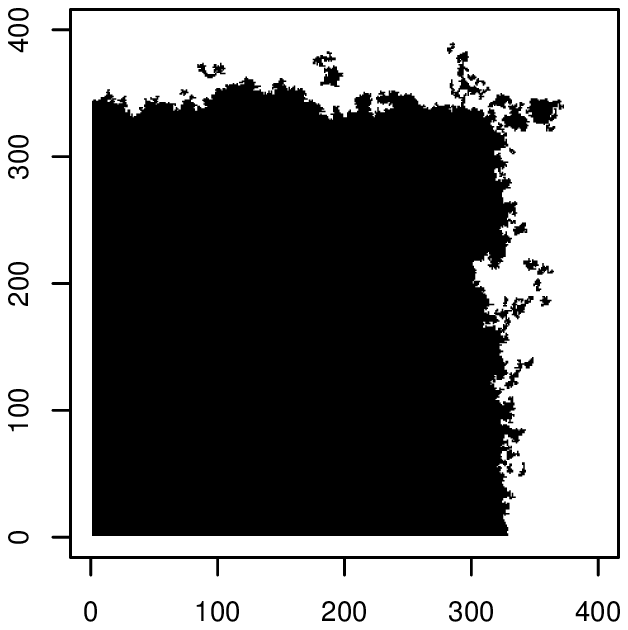}
\includegraphics[width=2.5in, height=2.5in]{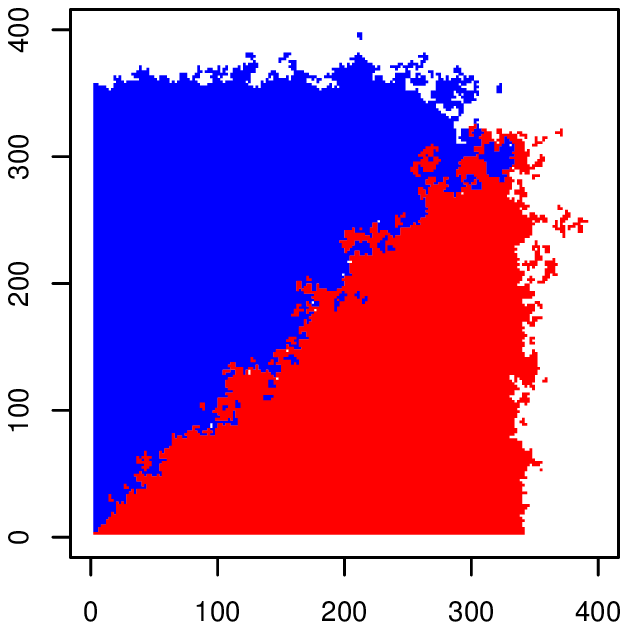} 
\end{tabular}
\caption{Growth and Competition Models in Hostile Environment.}
\end{figure}

%%%%%%%%%%%%%%%%%%%%%%%%%%%%%%%%%%%%%%%%%%%%%%%%%%%%%%%%%%%%%%%%%%%%%%%%%%%%%%%%%%%%%
Suppose that at time  zero every vertex of $\mathbb{Z}_{+}^2$ 
except the origin contains  a particle. 
There are two distinguished particles  located at $(1,0)$ and $(0,1)$, 
say Black particles.
All other vertices are occupied by  White particles.  
Every vertex flips  to the color of  a randomly chosen
south-west nearest neighbor with exponential rate $2$.
Thus, 
at time $t$ the color of a vertex $z$ is uniquely determined
by its  voter-admissible path. The set of Black particles $Q(t)$ 
is defined to be the set of all vertices $z$ such that the unique 
reverse voter-admissible path beginning at $(z,t)$ terminates at
$ \{ (1,0), (0,1) \} $.
Note that every vertex $(z_1,0)$ on the horizontal coordinate axes 
and every vertex $(0,z_2)$  on the vertical coordinate axes
eventually flips to Black color  and stays Black forever. 
Thus, almost surely for all large  $t$   vertex  $(z_1,z_2)$ 
is Black. 
By subadditivity, a shape theorem should hold for the 
growth model. 
Computer simulations of the growth model suggest that
the shape set is a square (see the first picture on  Figure \ref{Hostile}).
Below it is shown that  the limit shape  is exactly  $\mathcal{Q}$.
\begin{proposition} \label{growth:environment}
For every $\alpha \in (1/2,1)$ there exist  $c_1,c_2$ such that for all $t>0$
\[ P[ \mathcal{Q}(t-t^{\alpha}) \subset 
\hat{Q}(t)  
\subset \mathcal{Q}(t+t^{\alpha})]>
  1-c_1 t^2  \exp\{ -c_2 t^{(\alpha-1/2)} \}.  \]
\end{proposition}
\begin{proof}  
Recall that for every $t>0$ and $z \in \mathbb{Z}_{+}^2$, 
there exists a  unique reverse voter-admissible path $\tilde{ \gamma }_{(z,t)}$
 starting at $(z,t)$. The path travel downward, 
at rate $1$, and jumps across all inward-pointing arrows.  
Until the path hits the horizontal (vertical) axis  
the number of horizontal (vertical) jumps is
distributed as Poisson process with parameter $1$. 
Thus, there exist constants $c_1$ and $c_2$
such that  for every  $z \in Q(t-t^{\alpha})$
\[  
P( \tilde{\gamma}_{(z,t)}   \mbox{ terminates in } \{(1,0),(0,1) \} ) 
\ge 1- c_1 \exp \{ -c_2 t^{ (\alpha-1/2)}  \}.  
\]
For the same reason, there exist  constants $c_1$ and $c_2$  such that
for every $z \in Q^c(t+t^{\alpha})$,
\[  
P( \tilde{\gamma}_{(z,t)}   \mbox{ terminates in } \{ (1,0),(0,1) \} ) 
\le c_1 \exp \{ -c_2 t^{(\alpha-1/2)} \}.  
\]
The proposition follows from the fact that the  number of
vertices in  $Q(t-t^{\alpha})$ is of order at most $O(t^2)$
and the number of vertices on the boundary of 
 $Q(t+t^{\alpha})$ is of order at most $O(t)$.
\end{proof} 

If the growth models $Q(t)$ and $S(t)$ are coupled on the same
percolation structure $\Pi$, then clearly $Q(t) \subseteq S(t)$, and
thus $\mathcal{Q} \subseteq \mathcal{S}$.  Lemma
\ref{ProperSubset} asserts that $\mathcal{S}$ is strictly
larger than $\mathcal{Q}$.  

\begin{proof}
[Proof of Lemma \ref{ProperSubset}] The following argument was
communicated to the authors by Yuval Peres.  We consider a
representation of the Richardson model as a first passage percolation
model.  To each edge of the lattice associate a mean one exponential
random variable, also called a passage time of the edge.  The
variables are mutually independent.  For every pair of vertices
$z_1=(x_1,y_1)$, $z_2=(x_2,y_2)$ such that $x_1 \le x_2$ and $y_1 \le
y_2$ define the passage time $T(z_1,z_2)$ from $z_1$ to $z_2$ as the
infimum over traversal times of all North-East oriented paths from
$z_1$ to $z_2$.  The traversal time of an oriented path is the sum of
the passage times of its edges.  In the first passage percolation
description of the Richardson model, let
\[
Z(t)= \{ z \in \mathbb{Z}_{+}^2: 
T((1,0), z) \le t \mbox{  or }
T((0,1), z) \le t \}.
\]
It is enough to show that for some $\epsilon>0$, 
the vertex $(1,1)$ is in  $(1- \epsilon)\mathcal{S}$.
Consider a sequence of vertices  $z_n=(n,n)$  
on the main diagonal of the first quadrant  of $\mathbb{Z}^2$.
By the shape theorem, it  suffices to prove that almost surely
for infinitely many  $n$'s  the occupation times of $z_n$ satisfy
$T(z_n) \le n(1- \epsilon)$. \newline
\indent Consider vertices $(0,2),(2,0)$, and $(1,1)$. There are
exactly four oriented  distinct  paths  from the origin to these vertices.
Each such path has two edges  and expected passage time equal to $2$. 
Let $\gamma_{(1)}$ be the path with the smallest passage time among these four
paths. Denote by $X_1$ the terminal point of  $\gamma_{(1)}$,
 and denote its passage time  by $T_1$.
By symmetry $P(X_1=(0,2))= P(X_1=(2,0))= 1/4$ and $P(X_1=(1,2))=1/2$.
It easy to see that $ET_1<1$. Indeed,  let $\gamma_0$ be the 
path obtained by the following  procedure. Start at the origin
and make two oriented steps each time  moving in the direction
of the edge with minimal passage time (either north or east). 
Clearly 
\[
ET_1<E \tau(\gamma_0)=1
\]
where $\tau(\gamma_0)$ is the total passage time of $\gamma_0$. 
Restart at $X_1$ and repeat the procedure. Denote by $X_2$
the displacement on the second  step  
and by $T_2$ the passage time of the time minimizing path
from $X_1$  to $X_1+X_2$. 
Note that $W_k= \sum_{k=1}^{\infty} X_k$ is a random walk on 
$\mathbb{Z}_{+}^2$.
The random walk  visits  the main diagonal infinitely often in
such a way that $W_k=(k,k)$. 
Furthermore, if $S_k= \sum_{k=1}^{\infty} T_k $, 
then  by SLLN for some $\epsilon>0$
 almost surely  for all large $k$  we have $S_k \le (1-\epsilon)k$.
This finishes the proof.
\end{proof}
%%%%%%%%%%%%%%%%%%%%%%%%%%%%%%%%%%%%%%55
Suppose  now that at time zero there is one Red particle at 
$(1,0)$ and one Blue particle  at $(0,1)$.
As in the growth model, all other vertices are occupied
by White particles.
The flip rules are the same as in the growth model:
every vertex flips  to the color of a randomly chosen
south-west nearest neighbor with exponential rate $2$.
At time $t>0$ the color of a vertex $z$ is uniquely determined
by its  voter-admissible path. 
The  Red  cluster  $R^{(1)}(t)$ and the Blue cluster  $B^{(1)}(t)$ 
are  defined to be the sets of all vertices $z$ such that the unique 
reverse voter-admissible path beginning at $(z,t)$ terminates
respectively  at $(1,0)$  and $(0,1)$.
For $t>0$ define 
\[ K_1(t^{\alpha})= \{ z=(z_1,z_2)\in \mathbb{Z}_{+}^2 : z_1-z_2> t^{\alpha} \}, \] 
\[
K_2(t^{\alpha})= \{ z=(z_1,z_2)\in \mathbb{Z}_{+}^2 : z_2-z_1 > t^{\alpha} \}.
\]
\begin{proposition} \label{half-coloring:envir}
For every $\alpha \in (1/2,1)$ there exist  $c_1,\ c_2>0$ 
such that for all $t>0$
\begin{equation} \label{Triangle1}
 P[ \mathcal{Q}(t-t^{\alpha}) \cap  K_1(t^{\alpha})  
\subset \hat{R}^{(1)}(t) ] >
1-c_1 t^2  \exp\{ -c_2 t^{(\alpha-1/2)} \}, 
\end{equation}
\begin{equation} \label{Triangle2}
 P[ \mathcal{Q}(t-t^{\alpha}) \cap K_2(t^{\alpha})    
\subset \hat{B}^{(1)}(t) ] >
1-c_1 t^2  \exp\{ -c_2 t^{(\alpha-1/2)} \}.  
\end{equation}
\end{proposition}
\begin{proof}
We only show (\ref{Triangle1}), since  the proof of (\ref{Triangle2}) is identical.
First, observe that by Proposition \ref{growth:environment}
there exist $c_1, c_2>0$
\[ P[ \mathcal{Q}(t-t^{\alpha})  \subset  \hat{R}^{(1)}(t) \cup  \hat{B}^{(1)}(t)  ]
 >  1-c_1 t^2  \exp\{ -c_2 t^{(\alpha-1/2)} \}.  \]
Second, note that with  probability  exponentially close to one  
voter admissible paths of all  vertices
 $z \in  \mathcal{Q}(t-t^{\alpha}) \cap  K_1(t^{\alpha}) \cap \zz{Z}^2_{+} $
terminate below main diagonal.
That is, there exist  $c_1,c_2>0$ such that
\[  P \left( \tilde{\gamma}_{z,t)}(t)  \in  K_1(0)  \right)
  \ge 1-c_1 \exp \{-c_2 t^{(\alpha-1/2)}  \}.  \]
The result  (\ref{Triangle1})  immediately  follows from the above observations. 
\end{proof}
%%%%%%%%%%%%%%%%%%%%%%%%%%%%%%%%%%%%%%%%%%%%%%%%%%%%%%%%%%%%%%%%%%%%%%%%%%%%%%%%%%%%%%
%%%%%%%%%%%%%%%%%%%%%%%%%%%%%%%%%%%%%%%%%%%%%%%%%%%%%%%%%%%%%%%%%%%%%%%%%%%%%%%%%%%
\section{Oriented Competition model. Proof of Theorem \ref{MainTheorem}.}
 \label{sec:OrientedCompetition}
If the competition model
and  the competition model 
in hostile  environment  
are constructed on the same percolation structure $\Pi$, 
then almost surely for all $t>0$ 
\[  R^{(1)}(t) \subseteq R(t), \  B^{(1)}(t) \subseteq B(t). \]
Hence  it follows from Proposition \ref{half-coloring:envir} 
that almost surely for all large $t$
\[
\mathcal{Q}(t-t^{\alpha}) \cap  K_1(t^{\alpha}) \subset \hat{R}(t),
\]
\[
 \mathcal{Q}(t-t^{\alpha}) \cap K_2(t^{\alpha}) \subset \hat{B}(t). 
\]
Thus, asymptotically (as $t$ goes to infinity)  the  square 
$\mathcal{Q} \subset \mathcal{S}$ is colored deterministically. 
In particular, the region below the main diagonal is red, 
and the region above the diagonal is blue.
This proves the first part of  Theorem \ref{MainTheorem}. \newline
%%%%%%%%%%%%%%%%%%%%%%%%%%%%%%%%%%%%%%%%%%%%%%%%%%%%%%%%%%%%%%%%%%%%%%%%%
\indent 
The next question is what happens in the region 
$\mathcal{S} \setminus \mathcal{Q}$.
For each $z\in (\partial \mathcal{S})\setminus \mathcal{Q}$
and any $\varrho >0$, define the 
\emph{angular sector} $\mathcal{A} (z;\varrho) \subset \mathcal{S} \setminus \mathcal{Q}$
 of angular measure $\varrho$ rooted at $(1,1)$ and centered at $z$ by
\begin{gather*}
        \mathcal{A} (z;\varrho):= \{y \in \mathcal{S}
                \, : \, 
 |\mbox{arg} \{y-(1,1) \} - \mbox{arg} \{z-(1,1) \} | < \varrho/2  \}.
\end{gather*}
Fix $\epsilon>0$, $\alpha \in(1/2,1)$ and $\beta \in (1/2,1)$
such that $(\alpha +1)/2<\beta $.
For  $\varrho >0$ and $t \ge 1$,  let $A_{1}\subset A_{2}$
 be angular sectors with common center $z$ and angular measures
$r<r+t^{\beta-1}$, respectively, and such that
$ A_{2} \subset K_{\epsilon} $.
Define by $A_1^{c}$ and $A_2^{c}$ the complements of the sectors
in   $\mathcal{S} \setminus \mathcal{Q}$.
Fix $\delta \in (0,1)$, and set
\begin{align*}
        \mathcal{R}_{0}&=       \mathcal{R}_{0}^{t}=
         A_2(t-t^{\alpha}),\\
        \mathcal{B}_{0}&=       \mathcal{B}_{0}^{t}=
                  A_2^c (t+t^{\alpha}), \\
        \mathcal{B}_{1}&=       \mathcal{B}_{1}^{t}=                
                  A_1^c (t(1+\delta)+(t(1+\delta))^{\alpha}), \\
         \mathcal{R}_{1}&=       \mathcal{R}_{1}^{t}=
                 A_{1}(t(1+\delta)-(t(1+\delta))^{\alpha}).
\end{align*}

\begin{lemma}\label{lemma:stabilization}
There exist constants $c_{1},c_{2}>0$ such that the following is true,
for any $t\geq 1$. If the initial configuration $\xi ,\zeta$ is such
that $\hat{\xi} \supset \mathcal{R}_{0}^{t}$ and 
$\hat{\zeta} \subset \mathcal{B}_{0}^{t}$, then 
\begin{equation}\label{eq:stabilization}
        1-P_{\xi ,\zeta} [
                \hat{B} (\delta t)\subset \mathcal{B}_{1}^{t} ]
        \leq c_{1}t^{2} \exp \{-c_{2} (\delta  t)^{\alpha -1/2} \}.
\end{equation}
\end{lemma}
Lemma \ref{lemma:stabilization}
implies  that once an angular segment  is  occupied by 
one of the two types, it must remain so (except near its boundary) for 
a substantial amount of time.
Thus, Theorem \ref{MainTheorem} immediately follows from  Lemma  \ref{lemma:stabilization}  and  Theorem \ref{OrientedShapeTheorem}. 
For more details, see analogous  construction in \cite{kordzakh} (Section 4.3, pg. 14-15). \\
%%%%%%%%%%%%%%%%%%%%%%%%
\indent Let  $a=(1,1) \in \mathcal{S}$, 
be the right upper corner vertex of $\mathcal{Q}$,
and let  $b$ be a point on the boundary of  $\mathcal{S}$ 
such that $b \in  K_{\epsilon}$.
For $r,q  \in \mathbb{R}$ 
denote by  $I(r,q)$ an interval with ends at $r$ and $q$ 
and by $L(r,q)$ a line segment starting at  $r$ and passing through $q$.  
Let  $b'(\delta)$  be a point in the interval $I(a,b)$ 
 such that  $|b'-a| = |b-a| / (1+\delta) $  where $| \cdot |$ is an Euclidean norm.
For a point $r \in \mathbb{R}_{+}$  let $\overline{r}$ be the nearest vertex 
with integer coordinates. 
That is, $dist(r, \overline{r}) \le 1/2$  (if there is more than
one such vertex,  choose the vertex with the  smallest coordinates). \newline
\indent Suppose that at time zero the initial configuration
$R(0),B(0)$  is such that
$\hat{R}(0) \cup \hat{B}(0) \approx  \mathcal{S}t$  for large  $t>0$.
Then, by the shape theorem,    $\hat{R}(\delta t) \cup \hat{B}(\delta t) 
\approx  \mathcal{S}t(1+\delta)$.
Consider the  line segment $L(at(1+\delta),bt(1+\delta))$ starting at  
$at(1+\delta)$ and passing through $bt(1+\delta)$.  
Fix a point    $r \in L(a,b)$ such that  
$\overline{rt(1+\delta)}  \in  R(\delta t) \cup B(\delta t)$.
Note that  $rt(1+ \delta) \in L(at(1+\delta),bt(1+\delta))$.  
In  Claims  \ref{KeyClaim1},  \ref{KeyClaim2}, and  \ref{KeyClaim3} below,
it is shown that if  $\partial^{o}\mathcal{S}$  is uniformly curved,  
then with  probability exponentially close to one 
the ancestor of  $\overline{rt(1+\delta)}$ (if exists) is in 
the $(\delta t)^{\beta}$  neighborhood of  $I(at,bt)$  
for  some $\beta \in (3/4,1)$. 
In particular if $r \in I(b',b)$,
then the ancestor of  $\overline{rt(1+\delta)}$ (if exists) 
is in the $(\delta t)^{\beta}$  neighborhood of  $bt$.  
Observe that  this implies the statement of the  
Lemma  \ref{lemma:stabilization}. Consider three cases:
\begin{enumerate}
 \item \ $r \in I(a,b')$;
\item \ $r \in I(b', b)$ and $rt(1+\delta) \in  \mathcal{S}(t(1+\delta)-(t(1+\delta))^{\alpha})$; \item  \ $rt(1+\delta) \in  \mathcal{S}(t(1+\delta)+(t(1+\delta))^{\alpha}) 
\setminus  \mathcal{S}(t(1+\delta)-(t(1+\delta))^{\alpha})$. 
\end {enumerate}
The Claims  \ref{KeyClaim1},  \ref{KeyClaim2}, and  \ref{KeyClaim3} 
deal with the three cases respectively.
%%%%%%%%%%%%%%%%%%%%%%%%%%%%%%%%%%%%%%%%%%%%%%%%%%%%%%%%%%%%%%%%%%%%%%%%%%%
\begin{claim} \label{KeyClaim1}
% Fix  $\alpha \in (1/2,1)$.
There exist constants $c_1,c_2>0$  such that
for every  $b \in  \partial{ \mathcal{S}}  \cap K_{\epsilon}$ and for every 
 $r \in I(a,b')$,  
if the  initial  configuration   $R(0), B(0)$  is such that 
 $\mathcal{S}(t-t^{\alpha}) \subset   \hat{R}(0) \cup \hat{B}(0) \subset  \mathcal{S}(t+t^{\alpha})$   
then 
with probability at least  $1- c_1 \exp \{-c_2(\delta t)^{(\alpha- \frac{1}{2})} \}$ 
the ancestor of $\overline{r(1+\delta)t}$ exists and  is located in the  $(\delta t)^{\alpha}$
neighborhood of  $r(1+\delta)t- a\delta t \in I(at,bt)$.    
\end{claim}
\begin{proof}  
Recall that the voter admissible path is a continuous time random walk 
with exponential waiting times between jumps and drift $-a$.
By standard large deviations results, 
with probability exponentially close to one
the voter admissible reverse path  $\tilde{\gamma} \in \Gamma(\overline{r(1+\delta)t},\delta t)$ 
is attached to a vertex in     
the disk of radius $(\delta t)^{\alpha}$ centered at 
$(rt(1+\delta)- a \delta t) \in I(at,bt)$.
That is, for some  constants $c_1,\ c_2>0$, 
\[ P[ |  \tilde{\gamma}(t)  -  ( rt(1+\delta) - a \delta t)| > (\delta t)^{\alpha} ] 
\le c_1 \exp \{-c_2(\delta t)^{(\alpha- \frac{1}{2})}  \}.   \]
%Thus, the ancestor of $r$ is in the neighborhood of  $r'\in I(at,bt)$ 
\end{proof}
%%%%%%%%%%%%%%%%%%%%%%%%%%%%%%%%%%%%%%%%%%%%%%%%%%%%%%%%%%%%%%%%%%%%%%%%%%%%%%%%%%%
\begin{claim} \label{KeyClaim2}
% Fix $\alpha \in (1/2,1)$.
There exist constants $c_1,c_2>0$  such that
for every  $b  \in  \partial{ \mathcal{S} }  \cap K_{\epsilon}$
and for every $r \in I(b',b)$  with  
$rt(1+\delta) \in  \mathcal{S}(t(1+\delta)-(t(1+\delta))^{\alpha})$ 
if the  initial  configuration   $R(0), B(0)$  is such that 
$\mathcal{S}(t-t^{\alpha}) \subset   \hat{R}(0) \cup \hat{B}(0) \subset  \mathcal{S}(t+t^{\alpha})$,
then with probability at least  
$1- c_1 \exp \{-c_2(\delta t)^{(\alpha- \frac{1}{2})} \}$ 
the ancestor of $r(1+\delta)t$ exists and is located  
in the  $(\delta t)^{\beta}$ neighborhood of  $bt$.    
\end{claim}
\begin{proof}  
The heuristics of the proof are as follows. 
For  $t_1 \in (0, \delta t)$ consider 
a subset  $\Gamma_1 (\overline{rt(1+ \delta t)}, \delta t )$ 
of the set of reverse paths  $\Gamma (\overline{rt(1+ \delta)}, \delta t )$  
that contains only  those  paths that  coincide
with the reverse voter admissible path $\tilde{\gamma}$ 
on  $ \mathbb{Z}^2 \times (\delta t- t_1, \delta t)$. 
That is, for every  $\gamma  \in \Gamma_1 (\overline{rt(1+ \delta)}, \delta t )$,
for all $0<s<t_1$,    $\gamma(s) = \tilde{\gamma}(s)$ 
% (Clearly  $\tilde{\gamma} \in \Gamma_1 (rt(1+\delta), \delta t )$.)
The set of ends of  $\Gamma_1 (\overline{rt(1+ \delta)}, \delta t )$ is obtained 
by constructing reverse oriented Richardson process on 
 the subset $\mathbb{Z}_{+}^2 \times (0,t_2)$
of the percolation structure. The process starts with
 one occupied vertex at $\tilde{\gamma}(t_1)$, 
and runs  backward in time  for $t_2=\delta t-t_1$ units of time. 
By making an appropriate choice of  $t_1$  and $t_2$, we show that
with probability exponentially close to one   the  ancestor vertex of  
$(\overline{rt(1+\delta)}, \delta t)$ exists and is located  in the $(\delta t)^{\beta}$ 
neighborhood of  $bt$.
Denote by 
\[
\kappa = \frac{ |r-b| } { |b'-b| }.
\]
Consider  $ L( {\bf 0} , b(1+\delta)t )$, a line in $\mathbb{R}^2$  
connecting  the origin ${\bf 0}$ and  the point $b(1+ \delta )t$.
There exists a unique point  $r' \in   L({\bf 0}, b(1+\delta )t )$
between   $bt$ and   $b(1+ \delta )t$
such that 
\[
rt(1+ \delta) -r'= \left( \kappa \delta t \right) a
\]
\[
\mu(r'-bt)= (1- \kappa) \delta t.
\]
Set $t_1= \kappa \delta t - (\delta t)^{\alpha}$.
Note that on the percolation structure, if we start at 
$(\overline{rt(1+ \delta)}, \delta t)$, and  follow 
the reverse voter admissible  path for $t_1$ units of time, 
then with  probability exponentially close to one  
the end of the path is located in an Euclidean  disk with center  
at  $r'+  (\delta t)^{\alpha} a$
and radius $(\delta t)^{\alpha} \epsilon_1/4$
where $\epsilon_1$ is chosen so that
$a(1+\epsilon_1) \in \mathcal{S}$.
That is, for some  constants $c_1,\ c_2>0$, 
\[ P [ | \tilde{\gamma}(t_1) - (r'+  (\delta t)^{\alpha} a)| >  
(\delta t)^{\alpha} \epsilon_1/4 ] 
\le c_1 \exp \{-c_2(\delta t)^{(\alpha- \frac{1}{2})}  \}.   \]
Observe  also that if we  start a  reverse oriented Richardson process 
(i.e. South-West oriented Richardson process)
from any vertex $z$ in the   $(\delta t)^{\alpha} \epsilon_1 / 4 $ neighborhood 
of   $r'+  (\delta t)^{\alpha} a$, and run it backward  in time
for $t_2=(1-\kappa) \delta t+(\delta t)^{\alpha}$ units of time, then
\begin{equation}  \label{claim2:non-empty}
 P [ \Gamma ( z, t_2) \cap (R(0) \cup B(0)) = \emptyset ]
< c_1 \exp \{-c_2(\delta t)^{(\alpha- \frac{1}{2})}  \}.   
\end{equation}
Indeed, since 
\[
 \mathcal{S}(t- t^{\alpha}) \cap 
(z + \tilde{\mathcal{S}} (t_2- (\delta t)^{\alpha} \epsilon_1/4 )) \not=  \emptyset,
 \] 
(\ref{claim2:non-empty}) follows by   by Theorem~\ref{OrientedShapeTheorem}. \\
\indent Also, by Theorem  \ref{OrientedShapeTheorem} and  
by  Lemma \ref{OrientedShapeIntersection}, 
\[ P [ \Gamma ( z, t_2) \cap (R(0) \cup B(0))   
\not \subset D(bt, (\delta t)^{\beta}) ]
< c_1 \exp \{-c_2(\delta t)^{(\alpha- \frac{1}{2})}  \}.   \]
Thus,  with  probability exponentially close to one 
the intersection  of  the set 
 $\Gamma_1 (\overline{rt(1+ \delta t)}, \delta t )$ 
 with $R(0) \cup B(0)$
is non-empty and belongs to a disk 
of radius  $(\delta t)^{\beta}$
and center at $bt$.
 \end{proof}
\begin{claim} \label{KeyClaim3}
%  Fix $\alpha \in (1/2,1)$. 
There exist constants $c_1,c_2>0$  such that
for every  $b \in  \partial{ \mathcal{S}}  \cap K_{\epsilon}$ and 
for every  $r \in L(b',b)$  with 
$rt(1+\delta) \in  \mathcal{S}(t(1+\delta)+(t(1+\delta))^{\alpha}) 
\setminus  \mathcal{S}(t(1+\delta)-(t(1+\delta))^{\alpha})$, 
if the  initial  configuration   $R(0),B(0)$  is such that 
 $\mathcal{S}(t-t^{\alpha}) \subset  
\hat{R}(0)  \cup \hat{B}(0)  
\subset  \mathcal{S}(t+t^{\alpha})$,
then  with probability at least  
$1- c_1 \exp \{-c_2(\delta t)^{(\alpha- \frac{1}{2})} \}$ 
the set of potential ancestors of $\overline{r(1+\delta)t}$ 
is either empty or it is contained 
in the   $(\delta t)^{\beta}$ neighborhood of  $bt$.    
\end{claim}
\begin{proof}
If the ancestor of $\overline{rt(1+ \delta)}$ exists, 
it is located in the set of ends of  $\Gamma(\overline{rt(1+ \delta)}, \delta t )$.
The set of ends of  $\Gamma(\overline{rt(1+ \delta)}, \delta t )$ 
is obtained by constructing reverse oriented Richardson process 
starting with one occupied
vertex at $\overline{rt(1+ \delta)}$, and running the process  
on  the subset $\mathbb{Z}_{+}^2 \times (0,\delta t)$  of 
the percolation structure  
backward in time  for $\delta t$ units of time. 
Then by Theorem~\ref{OrientedShapeTheorem}  and 
Lemma \ref{OrientedShapeIntersection},  
\[ P [  \Gamma( \overline{rt(1+ \delta)}, \delta t) \cap (R(0) \cup B(0))  
  \not \subset D(bt, (\delta t)^{\beta}) ]
\le c_1 \exp \{-c_2(\delta t)^{(\alpha- \frac{1}{2})}  \}.   \]
\end{proof}
The Claims \ref{KeyClaim1}, \ref{KeyClaim2} and \ref{KeyClaim3},
imply  the statement of Lemma \ref{lemma:stabilization}.
This  finishes the proof of  Theorem  \ref{MainTheorem}.
%%%%%%%%%%%%%%%%%%%%%%%%%%%%%%%%%%%%%%%%%%%%%%%%%%%%%%%%%%%%%%%%%%%%%%
%%%%%%%%%%%%%%%%%%%%%%%%%%%%%%%%%%%%%%%%%%%%%%%%%%%%%%%%%%%%%%%%%%%%%%
\section*{Acknowledgment}  
We thank Yuval Peres for helpful conversations and for
suggesting the  proof of Lemma \ref{ProperSubset}.

\end{document}